\documentclass{elsarticle}

\usepackage{amsmath}
\usepackage{amsthm}
\usepackage{amssymb}
\usepackage{upref}
\usepackage{amscd}

\numberwithin{equation}{section}

\newtheorem{theorem}[equation]{Theorem}
\newtheorem{lemma}[equation]{Lemma}
\newtheorem{proposition}[equation]{Proposition}
\newtheorem{corollary}[equation]{Corollary}

\theoremstyle{definition}

\DeclareMathOperator{\ann}{ann}

\DeclareMathOperator{\Hom}{Hom}

\DeclareMathOperator{\End}{End}

\DeclareMathOperator{\bimod}{Bimod}

\newcommand{\Q}{\mathbb{Q}}
\newcommand{\ideal}[1]{\mathfrak{#1}}

\newcommand{\cat}[1]{\mathcal{#1}}

\newcommand{\Z}{\mathbb{Z}}

\newcommand{\Mod}{\text{-mod}}

\newcommand{\mathcolon}{\colon\,}

\newcommand{\uc}{\textup{:}}

\begin{document}
 
\title{Additive closed symmetric monoidal structures on $R$-modules}

\date{\today}

\author{Mark Hovey}
\address{Department of Mathematics \\ Wesleyan University
\\ Middletown, CT 06459}
\ead{hovey@member.ams.org}

\begin{keyword}
symmetric monoidal, closed symmetric monoidal, module
\MSC 18D10 \sep 16D90 
\end{keyword}

\begin{abstract}
In this paper, we classify additive closed symmetric monoidal
structures on the category of left $R$-modules by using Watts'
theorem.  An additive closed symmetric monoidal structure is
equivalent to an $R$-module $\Lambda _{A,B}$ equipped with two
commuting right $R$-module structures represented by the symbols $A$
and $B$, an $R$-module $K$ to serve as the unit, and certain
isomorphisms.  We use this result to look at simple cases.  We find
rings $R$ for which there are no additive closed symmetric monoidal
structures on $R$-modules, for which there is exactly one (up to
isomorphism), for which there are exactly seven, and for which there
are a proper class of isomorphism classes of such structures.  We also
prove some general structual results; for example, we prove that the
unit $K$ must always be a finitely generated $R$-module.
\end{abstract}

\maketitle

\section*{Introduction}

It is well known that the category of left $R$-modules becomes closed
symmetric monoidal under the tensor product $A\otimes_{R}B$ if and
only if $R$ is commutative.  However, there are many other cases when
the category of $R$-modules is closed symmetric monoidal.  For
example, if $k$ is a field and $G$ is a group, the category of
$k[G]$-modules (that is, representations of $G$ on $k$-vector spaces)
is closed symmetric monoidal under $A\otimes_{k}B$, even though $k[G]$
is not commutative in general.  This is explained by the fact that
$k[G]$ is a Hopf algebra.  But there are other examples where $R$ is
not a Hopf algebra, such as the category of perverse $R$-modules
considered in~\cite{hovey-intersection}.  

So a natural question to ask is just what one needs to know about $R$
in order to produce a closed symmetric monoidal structure on the
category of $R$-modules.  Of course, we do not really want an
arbitrary closed symmetric monoidal structure; we require that the
monoidal product be an additive functor in both variables.  We would
like to be able to answer basic questions such as the following.  Are
there rings $R$ where the category of $R$-modules cannot be given an
additive closed symmetric monoidal structure?  Are there rings $R$
where the category of $R$-modules possesses a unique additive closed
symmetric monoidal structure?

At first glance, such problems seem completely intractable because
closed symmetric monoidal structures are so complicated, involving the
entire category of $R$-modules.  The key ingredient, though, is Watts'
theorem~\cite{watts}.  This theorem says that any additive functor $F$
from $R$-modules to abelian groups that is right exact and commutes
with direct sums is naturally isomorphic to $\Lambda \otimes_{R} (-)$
for some $R$-bimodule $\Lambda$.  After some work, we then see that a
closed symmetric monoidal structure $A\wedge B$ on left $R$-modules
must be given by
\[
A\wedge B\cong ((R\wedge R)\otimes_{R}A)\otimes_{R}B,
\]
so that the functor $-\wedge -$ is determined by $R\wedge R$, as a
$2$-fold bimodule (one left module structure, and two right module
structures).  

Then the natural thing to do is try to determine which $2$-fold
bimodules $\Lambda_{A,B}$ actually arise as $R\wedge R$ for some
closed symmetric monoidal structure $-\wedge -$.  This is more
complicated than it seems because one must deal with the coherence
isomorphisms of a closed symmetric monoidal structure, but of course
it can be done.  We also determine when two symmetric monoidal
structures determined by $\Lambda$ and $\Gamma$, respectively, are
equivalent as symmetric monoidal functors.  This involves an
isomorphism 
\[
\Gamma \otimes X\otimes X\xrightarrow{}X\otimes \Lambda 
\]
of $2$-fold bimodules, where $X$ is an element of the bimodule Picard
group of $R$.  

We establish some basic structural results, though we think there is
much more to say.  For example, we show that the unit $K$ of an
additive closed symmetric monoidal structure on left $R$-modules must
be a finitely generated $R$-module with a commutative endomorphism
ring.  To proceed further along these lines, it might be worthwhile to
develop a theory of flatness for an additive symmetric monoidal
structure $\wedge$ on $R$-modules, and concentrate on those additive
symmetric monoidal structures for which projectives are flat.  

We also consider examples.  For example, if $R$ is a field or a
principal ideal domain that does not contain a field, then there is
exactly one additive closed symmetric monoidal structure on
$R$-modules (up to symmetric monoidal equivalence).  If $R$ is a
division ring that is not a field, there are no additive closed
symmetric monoidal structures on $R$-modules.  If $R$ is the group
ring $\mathbb{F}_{2}[\mathbb{Z}/2]$, there are precisely $7$ different
closed symmetric monoidal structures on $R$-modules, though only three
different underlying functors.  If $R$ is the group ring
$k[\mathbb{Z}/2]$ where the characteristic is not $2$, however, there
are a proper class of inequivalent additive closed symmetric monoidal
structures on the category of $R$-modules.  Most of these cannot come
from Hopf algebra structures on $R$.  

Throughout this paper, the symbol $\otimes$ will denote the tensor
product over the ring $R$ unless otherwise stated.  Furthermore, all
functors will be assumed to be additive, even if not explicitly stated
to be so.

\section{$n$-fold bimodules}\label{sec-notation}

Throughout this paper, we will be working with the category of left
$R$-modules, but we will frequently have to work with left $R$-modules
that have multiple different commuting right $R$-module structures.
This necessitates some complicated notation.  We will denote an
$R\otimes_{\Z } (R^{\textup{op}})^{\otimes_{\Z } n}$-module by
$\Lambda_{1,2,\dotsc ,n}$, where the subscripts denote the commuting
right $R$-actions.  If we need elements, we will write $x\odot_{i}r$
for the $i$th multiplication.  We will call such an object an $n$-fold
bimodule, and denote the category of such things as $\bimod_{n} (R)$.
Now the symmetric group $\Sigma_{n}$ acts on $\bimod_{n} (R)$ by
permuting the right module structures.  Indeed, if $\sigma \in
\Sigma_{n}$ is a permutation, then
\[
(\sigma \Lambda)_{1,2,\dotsc ,n} = \Lambda_{\sigma (1),\dotsc ,\sigma
(n)}.  
\]
We note that in practice, as we will see below, it is usually easier
to denote the different right module structures by letters, such as
$\Lambda_{A,B,C}$.  If $\sigma$ denotes the permutation $(1 3 2)$, for
example, then we would write $\sigma \Lambda$ as $\Lambda_{C,A,B}$.  A
map
\[
f\mathcolon \Lambda_{A,B,C}\xrightarrow{}\Lambda_{C,A,B}
\]
would then be a map of left $R$-modules such that $f (x\odot_{1}r)=f
(x)\odot_{2}r$ (matching up the position of the $A$'s), $f
(x\odot_{2}r)=f (x)\odot_{3}r$, and $f (x\odot_{3}r)=f
(x)\odot_{1}r$.  

Now, note that the tensor product over $R$ defines a bifunctor
\[
\bimod_{n} (R)\times \bimod_{m} (R)\xrightarrow{}\bimod_{n+m-1} (R),
\]
obtained by tensoring the $n$th right module structure on the first
factor with the left module structure on the second factor.  For
example, we will have expressions like 
\[
f\mathcolon \Lambda_{C,\Gamma } \otimes \Gamma _{B,A} \xrightarrow{}
\Lambda_{\Gamma ,A}\otimes \Gamma _{C,B}.   
\]
This denotes a map of $3$-fold bimodules.  In the domain, we tensor
the second right $R$-module structure on $\Lambda$ with the left
module structure on $\Gamma$, whereas in the range we tensor the first
bimodule structure on $\Lambda$ with the left module structure on
$\Gamma$.  Furthermore, we have $f (( x\odot_{1}r) \otimes y)=f (x\otimes
y)\odot_{2}r$, where the second right module structure on the range comes
from the first right module structure on $\Gamma $.  Similarly, we
have $f (x\otimes (y\odot_{1}r))=f (x\otimes y)\odot_{3}r$.  We also have $f
(x\otimes (y\odot_{2}r))=f (x\otimes y)\odot_{1}r$, where the first
right module structure on the range comes from the \textbf{second}
right module structure on $\Lambda$, because we have used the first
one to form the tensor product in the range.

\section{Closed symmetric monoidal structures}\label{sec-main}

In this section, we prove our main classification result.  We remind
the reader that all tensor products are over $R$ unless otherwise
stated.  

\begin{theorem}\label{thm-watts}
Suppose the category of $R$-modules admits an additive closed
symmetric monoidal structure $-\wedge -$.  Then $R\wedge R$ is a
$2$-fold bimodule, and there is a natural isomorphism of bifunctors
\[
(R\wedge R)_{B,A}\otimes A\otimes  B \cong A\wedge B.
\]
\end{theorem}

This theorem is written using the notation of the previous section, so
that, in the domain of this isomorphism 
\[
(x\odot_{1}r)\otimes y\otimes z= x\otimes y\otimes (rz) \text{ and }
(x\odot_{2}r)\otimes y\otimes z = x\otimes (ry)\otimes z.
\]
Also note that if $R$ is a $k$-algebra, for a commutative ring $k$, we
can look at $k$-linear closed symmetric monoidal structures.  This
would mean that multiplication by $x\in k$ on $A$ would induce
multiplication by $x$ on $A\wedge B$ and on $B\wedge A$ for all $B$.
In this case, the three potentially different actions of $k$ on
$R\wedge R$ would in fact all be the same.  

\begin{proof}
Fix an $R$-module $A$, and consider the functor $B\mapsto A\wedge B$
from left $R$-modules to left $R$-modules.  This functor preserves
direct sums and is right exact (it is a left adjoint, because of the
closed structure).  Watt's theorem~\cite{watts} then implies that
$A\wedge R$ is an $R$-bimodule, and 
\[
(A\wedge R)\otimes B\cong A\wedge B
\]
naturally in $B$.  To see that this is also natural in $A$, recall
that the map 
\[
\alpha_{A,B}\mathcolon (A\wedge R)\otimes B\xrightarrow{}A\wedge B
\]
is defined as follows, in the proof of Watt's theorem.  Given $b\in
B$, let $\phi_{b}\mathcolon R\xrightarrow{}B$ denote the map of
$R$-modules that takes $1$ to $b$.  There is then an induced map
$A\wedge \phi_{b}$, and $\alpha_{A,B} (x\otimes b)$ is defined to be
$(A\wedge \phi_{b}) (x)$.  From this, it is easy to check that
$\alpha_{A,B}$ is natural in $A$ as well.

Now let $G (A)=A\wedge R$.  Then $G (A)$ is a right exact functor from
$R$-modules to $R$-bimodules (or
$R\otimes_{\Z}R^{\textup{op}}$-modules) that preserves direct sums.
We can therefore apply Watt's theorem~\cite{watts} again to give us
the desired result.
\end{proof}

The natural question then arises as to which $2$-fold bimodules
$\Lambda_{B,A}$ define a closed symmetric monoidal structure on the
category of $R$-modules.

We first point out that the closed structure always exists.  

\begin{lemma}\label{lem-closed}
Suppose $\Lambda_{B,A}$ is a $2$-fold $R$-bimodule.  Then we have a
natural isomorphism
\[
R\Mod (\Lambda_{B,A}\otimes A\otimes B, P) \cong R\Mod (A,
\bimod (\Lambda_{B,A},\Z \Mod (B,P))),
\]
where we use the right $R$-module structure on $\Lambda $ denoted by
the subscript $B$ to form $\bimod (\Lambda_{B,A},\Z \Mod (B,P))$, and
we use the one denoted by the subscript $A$ to make this abelian group
into a left $R$-module.
\end{lemma}

Here $\Hom_{\Z }(B,P)$ is an $R$-module via the $R$-action on $P$, and
a right $R$-module via the $R$-action on $B$.

\begin{proof}
This is really just an exercise in adjointness of tensor and Hom,
though one has to be careful to keep track of all the actions.  It is
easiest to work more generally.  Suppose $M$ is a bimodule.  Then we
have a natural isomorphism
\[
\phi \mathcolon R\Mod (M\otimes B,P)\cong \bimod (M, \Z \Mod (B,P)).
\]
This isomorphism is defined as usual by $\phi (f) (m) (b)=f (m\otimes
b)$.  The reader must check that $\phi (f)$ is a map of bimodules.  To
see that $\phi$ is an isomorphism, one constructs its inverse $\psi$,
where $\psi (g) (m\otimes b)=g (m) (b)$.  Again, there are many
details to check, which we leave to the reader.  Applying this
isomorphism to $A=\Lambda_{B,A}\otimes A$, we get
\[
R\Mod (\Lambda_{B,A}\otimes A\otimes B,P)\cong \bimod
(\Lambda_{B,A}\otimes A,\Z \Mod (B,P)).
\]

Now suppose $N$ is a general bimodule.  Then there is a natural
isomorphism
\[
\sigma \mathcolon \bimod (\Lambda_{B,A}\otimes A,N)
\xrightarrow{} R\Mod (A,\bimod (\Lambda_{B,A},N)).  
\]
Once again, we have $\sigma (f) (a) (\lambda )=f (\lambda \otimes a)$,
and, for the inverse $\tau$ of $\sigma$, we have $\tau (g) (\lambda \otimes
a)=g (a) (\lambda )$.  We leave it to the reader to check the details.
Taking $N=\Hom (B,P)$ completes the proof.
\end{proof}

Naturally, the other conditions necessary for a symmetric monoidal
structure are considerably more complicated.  The basic idea, however, is
simple.  In order to get symmetry of the product $A\wedge B$, we need
the two right module structures on $\Lambda $ to be isomorphic.
In order to get associativity of $A\wedge B$, we need the three
different right module structures on $\Lambda \otimes \Lambda$ to be
isomorphic.  These two together, of course, will imply that all of the
different permutations of the $n+1$ different right module structures
on $\Lambda^{\otimes n}$ will be isomorphic.  Then we also need a
unit.  

\begin{theorem}\label{thm-coherence}
Let $\Lambda_{B,A}$ be a $2$-fold bimodule used to define $-\wedge -$
on $R$-modules.  There is a one-to-one correspondence between additive
closed symmetric monoidal structures on $R$-modules with $-\wedge -$
as the monoidal product and the following data\uc
\begin{enumerate}
\item [(a)] An associativity isomorphism 
\[
a\mathcolon \Lambda_{C,\Lambda }\otimes \Lambda_{B,A}\xrightarrow{}
\Lambda_{\Lambda ,A}\otimes \Lambda_{C,B}.
\]
This can be remembered by noting that the subscripts on the first
$\Lambda$ in the target are the second subscripts on the two
$\Lambda$'s in the domain, and the subscripts on the second $\Lambda$
in the target are the first subscripts on the two $\Lambda$'s in the
domain.  
\item [(b)] A left $R$-module $K$ and a unit isomorphism
$\ell \mathcolon \Lambda_{B,K}\otimes K\cong R_{B}$ of bimodules.  
\item [(c)] A commutativity isomorphism $c\mathcolon
\Lambda_{B,A}\xrightarrow{}\Lambda_{A,B}$.
\end{enumerate}
This data must satisfy the following coherence conditions.  
\begin{enumerate}
\item (Associativity pentagon) Let $\Gamma =\Delta =\Lambda $ for
notational clarity.  Then the
following diagram commutes.
\small
\[
\minCDarrowwidth20pt
\begin{CD}
\Lambda _{D,\Gamma }\otimes \Gamma _{C,\Delta } \otimes \Delta _{B,A}
@>a\otimes 1>> \Lambda _{\Gamma ,\Delta } \otimes \Gamma _{D,C}
\otimes \Delta _{B,A} @>1\otimes T>>
\Lambda_{\Delta ,\Gamma } \otimes \Gamma_{B,A} \otimes \Delta_{D,C} \\
@V 1\otimes a VV @. @VVa\otimes 1 V \\
\Lambda_{D,\Gamma} \otimes \Gamma_{\Delta ,A} \otimes \Delta_{C,B}
@>>a\otimes 1> \Lambda_{\Gamma ,A} \otimes \Gamma_{D, \Delta} \otimes
\Delta_{C,B} @>>1\otimes a> \Lambda_{\Gamma ,A} \otimes \Gamma_{\Delta
,B} \otimes \Delta_{D,C}
\end{CD}
\]
Here $1\otimes T$ switches the last two factors using the
commutativity isomorphism of $\otimes_{R}$, but also reverses the
symbols $\Gamma$ and $\Delta$, which after all both mean $\Lambda$.
This necessitates changing the subscripts on $\Lambda$ as well.  
\item (Compatibility of left and right unit) The following diagram
commutes:
\[
\begin{CD}
\Lambda _{B,\Lambda } \otimes \Lambda _{K,A}\otimes K @>a\otimes 1>>
\Lambda _{\Lambda ,A}\otimes \Lambda _{B,K}
\otimes K @>1\otimes \ell >> \Lambda _{B,A} \\
@| @. @| \\
\Lambda _{B,\Lambda } \otimes \Lambda _{K,A}\otimes K @>>1\otimes c\otimes 1>
\Lambda _{B,\Lambda }\otimes \Lambda _{A,K} \otimes K @>>1\otimes \ell > \Lambda _{B,A}.
\end{CD}
\]
\item (Commutativity-associativity hexagon) The following diagram commutes. 
\small
\[
\minCDarrowwidth20pt
\begin{CD}
\Lambda _{C,\Lambda }\otimes \Lambda_{B,A} @>c>>
\Lambda_{C,\Lambda} \otimes \Lambda_{A,B} @>a>>
\Lambda_{\Lambda,B}\otimes \Lambda_{C,A} @>c>> \Lambda_{\Lambda,B}\otimes \Lambda_{A,C} \\
@| @. @. @| \\
\Lambda_{C,\Lambda}\otimes \Lambda_{B,A} @>>a>
\Lambda_{\Lambda,A}\otimes \Lambda_{C,B} @>>c>
\Lambda_{A,\Lambda}\otimes \Lambda_{C,B} @>>a> \Lambda_{\Lambda,B}
\otimes \Lambda_{A,C}.
\end{CD}
\]
\item The composite 
\[
\Lambda_{B,A} \xrightarrow{c} \Lambda_{A,B} \xrightarrow{c} \Lambda_{B,A}
\]
is the identity.  
\end{enumerate}
\end{theorem}

If $R$ is a $k$-algebra for a commutative ring $k$, and we are looking
at $k$-linear closed symmetric monoidal structures, then this theorem
remains true as long as the three different $k$-module structures on
$\Lambda$ are the same.  

\begin{proof}
Just using the usual associativity and commutativity isomorphisms for
$\otimes_{R}$, we find natural isomorphisms 
\[
(A\wedge B)\wedge C\cong \Lambda_{C.\Lambda} \otimes \Lambda_{B,A}
\otimes A\otimes B\otimes C 
\]
and 
\[
A\wedge (B\wedge C) \cong \Lambda_{\Lambda ,A} \otimes \Lambda_{C,B}
\otimes A \otimes B \otimes C.  
\]
Given $a$, it is now clear how to define a natural associativity
isomorphism $a_{A,B,C}$ for $-\wedge -$, simply as $a\otimes 1\otimes
1\otimes 1$.  On the other hand, given the natural associativity
isomorphism $a_{A,B,C}$, we let $A=B=C=R$ to get $a$.  One can see
that $a$ then respects the given right module structures by using
naturality with respect to the right multiplication by $x$ maps
$r_{x}\mathcolon R\xrightarrow{}R$, in the $A$, $B$, and $C$ slots.  

There is a similar equivalence between the isomorphism $\ell
\mathcolon \Lambda_{B,K}\otimes K\xrightarrow{}R_{B}$ and a natural
left unit isomorphism $\ell_{B}\mathcolon K\wedge B \xrightarrow{}B$.
There is also a similar equivalence between $c$ and a natural
commutativity isomorphism $c_{A,B}$.

An excellent reference for the coherence diagrams needed to make
$a_{A,B,C}$, $\ell_{B}$, and $c_{A,B}$ part of a symmetric monoidal
structure is~\cite{joyal-street}, particularly Propositions~1.1
and~2.1.  They show that the only coherence diagrams needed are the
associativity pentagon, the compatibility between the right and left
unit, the commutativity-associativity hexagon, the fact that the right
unit is $r_{B}=\ell_{B}c_{B,K}$, and the fact that $c^{2}$ is the
identity.  Given $c_{A,B}$, this means we do not need $r_{B}$, so we
have omitted it.  One must now merely translate these coherence
diagrams into analogous facts about $a$, $\ell$, and $c$ to complete
the proof.  

The associativity pentagon is perhaps the most confusing, so we will
discuss that one in some detail, and leave the others to the reader.
Here is the standard associativity pentagon.  
\small
\[
\begin{CD}
((A\wedge B)\wedge C)\wedge D @>a_{A\wedge B,C,D}>> (A\wedge B)\wedge
(C\wedge D) @>a_{A,B,C\wedge D}>> A\wedge (B\wedge (C\wedge D)) \\
@V a_{A,B,C}\wedge 1 VV @. @| \\
(A\wedge (B\wedge C))\wedge D @>> a_{A,B\wedge C,D} > A\wedge
((B\wedge C)\wedge D) @>> 1\wedge a_{B,C,D} > A\wedge (B\wedge (C\wedge D))
\end{CD}
\]
\normalsize
Using the standard commutativity and associativity isomorphisms of
$\otimes$, the first term 
\[
((A\wedge B)\wedge C)\wedge D
\]
is represented by $\Lambda_{D,\Gamma} \otimes \Gamma_{C,\Delta}\otimes
\Delta_{B,A}$ (tensored with $A\otimes B\otimes C\otimes D$).  Here
$\Delta_{B,A}$ tensors $A$ and $B$, $\Gamma_{C,\Delta}$ tensors
$(A\wedge B)$ and $C$, and $\Lambda_{D, \Gamma}$ tensors $((A\wedge
B)\wedge C)$ and $D$.  The map $a_{A\wedge B,C,D}$ treats $A\wedge B$
as a single object, and will therefore leave the factor $\Delta_{B,A}$
unchanged, so is represented by 
\[
a\otimes 1 \mathcolon \Lambda _{D,\Gamma }\otimes \Gamma _{C,\Delta }
\otimes \Delta _{B,A} \xrightarrow{} \Lambda _{\Gamma ,\Delta }
\otimes \Gamma _{D,C} \otimes \Delta _{B,A}.
\]
The map $a_{A,B, C\wedge D}$ treats $C\wedge D$ as a single object,
and will therefore leave $\Gamma_{D,C}$ unchanged and apply $a$ to the
other two.  In order to do this, we would like to switch the order of
$\Gamma_{D,C}$ and $\Delta_{B,A}$, and then apply $a\otimes 1$.  It
turns out to be notationally much easier later if we also reverse the
names of $\Gamma$ and $\Delta$, so that the map $a_{A,B,C\wedge D}$ is
represented by the composite 
\[
\Lambda _{\Gamma ,\Delta }
\otimes \Gamma _{D,C} \otimes \Delta _{B,A} \xrightarrow{1\otimes T}
\Lambda_{\Delta ,\Gamma} \otimes \Gamma_{B,A} \otimes \Delta_{D,C}
\xrightarrow{a\otimes 1} \Lambda_{\Gamma ,A} \otimes \Gamma_{\Delta
,B} \otimes \Delta_{D,C}.  
\]  
This completes the clockwise half of the associativity pentagon. The
counterclockwise part is simpler. The first map $a_{A,B,C}\wedge 1$
leaves $D$ alone, so will leave $\Lambda_{D,\Gamma}$ alone.  It is
therefore represented by 
\[
1\otimes a \mathcolon \Lambda_{D, \Gamma} \otimes \Gamma_{C,\Delta}
\otimes \Delta _{B,A} \xrightarrow{} \Lambda_{D,\Gamma} \otimes
\Gamma_{\Delta ,A} \otimes \Delta _{C,B}.
\]
The next map $a_{A,B\wedge C,D}$ treats $B\wedge C$ as a single
entity, so will leave $\Delta _{C,B}$ alone.  It is then represented
by 
\[
a\otimes 1 \mathcolon \Lambda_{D,\Gamma} \otimes
\Gamma_{\Delta ,A} \otimes \Delta _{C,B}\xrightarrow{} \Lambda_{\Gamma ,A} \otimes
\Gamma_{D, \Delta } \otimes \Delta _{C,B}.
\]
Finally, the last map $1\otimes a_{B,C,D}$ leaves $A$ alone, so will
leave $\Lambda_{\Gamma ,A}$ alone.  It is reprensented by 
\[
a\otimes 1\mathcolon \Lambda_{\Gamma ,A} \otimes
\Gamma_{D, \Delta } \otimes \Delta _{C,B} \xrightarrow{} \Lambda_{\Gamma ,A} \otimes
\Gamma_{\Delta ,B} \otimes \Delta _{D,C}.  
\]
This completes the construction of the associativity pentagon.  
\end{proof}

We now point out that Watt's theorem can also be used to classify
additive symmetric monoidal equivalences between additive symmetric
monoidal structures on $R$-modules.  In an attempt to make the various
bimodule structures clear, we have used $Y$ and $Z$ as alternative
names for $X$ in the theorem below.  We have also used $T$ for the
usual commutativity isomorphism of the tensor product and for a
general permutation of tensor factors.

\begin{theorem}\label{thm-equiv}
Suppose $\wedge $ and $\Box$ are additive symmetric monoidal
structures on the category of $R$-modules with units $K$ and $K'$,
respectively, and represented by the $2$-fold bimodules $\Lambda $ and
$\Gamma $, respectively.  Then an additive symmetric monoidal functor
from $\wedge $ to $\Box $ that has a right adjoint is equivalent to an
bimodule $X$, an isomorphism $\eta \mathcolon K'\xrightarrow{}X\otimes
K$, and an isomorphism
\[
m\mathcolon \Gamma_{Y,X}\otimes X\otimes Y \xrightarrow{} X\otimes
\Lambda 
\]
of $2$-fold bimodules, such that the following diagrams commute:
\begin{enumerate}
\item (Unit) 
\small
\[
\begin{CD}
\Gamma_{Y,K'} \otimes K' \otimes Y @>1\otimes \eta \otimes 1 >>
\Gamma_{Y,X} \otimes X_{K}\otimes K \otimes Y @>1\otimes 1\otimes T>> \Gamma_{Y,X} \otimes
X_{K}\otimes Y\otimes K \\
@V\ell \otimes 1 VV @. @VVmV \\
X @= X @<<1\otimes \ell <  X\otimes \Lambda_{Y,K} \otimes K 
\end{CD}
\]
\item (Commutativity) 
\[
\begin{CD}
\Gamma_{Y,X} \otimes X\otimes Y @>m>> X \otimes \Lambda \\
@Vc\otimes T VV @VV1\otimes c V \\
\Gamma_{X,Y} \otimes Y\otimes X @>>m> X\otimes \Lambda 
\end{CD}
\]
\item (Associativity)
\small
\[
\minCDarrowwidth15pt
\begin{CD} 
\Gamma_{Z,\Gamma }\otimes \Gamma_{Y,X} \otimes X\otimes  Y \otimes Z
@>m>> \Gamma_{Z,X} \otimes (X\otimes \Lambda) \otimes Z
@>T>>
\Gamma_{Z,X} \otimes X\otimes Z\otimes \Lambda \\ 
@Va\otimes T VV @. @VVm\otimes 1 V \\
\Gamma_{\Gamma ,X} \otimes \Gamma_{Z,Y} \otimes Y \otimes Z \otimes X @. \
@. X\otimes \Lambda \otimes \Lambda  \\
@V1 \otimes m\otimes 1 VV @. @VV 1\otimes a V \\
\Gamma_{X,X} \otimes (X\otimes \Lambda) \otimes X @>>T> \Gamma_{X,X}
\otimes X\otimes X\otimes \Lambda @>>m> X \otimes \Lambda \otimes
\Lambda
\end{CD}
\]
\end{enumerate}
\end{theorem}

Composition of additive symmetric monoidal functors corresponds to the
tensor product of bimodules, and the identity functor corresponds to
the bimodule $R_{R}$.  Thus, additive symmetric monoidal equivalences
of additive symmetric monoidal structures are given by tensoring with
a bimodule that lies in the bimodule Picard group
(see~\cite{yekutieli}). In fact, if $X$ lies in the bimodule Picard
group, then tensoring with $X$ loses no information.  In this case,
then, the compatibility diagrams above show that the isomorphisms
$\ell ,c$, and $a$ for $\Box$ are determined by the corresponding
isomorphisms for $\wedge$, $m$, and $\eta$.  Thus we can think of the
bimodule Picard group as acting on symmetric monoidal structures with
fixed unit $K$, though there is also an action by the automorphisms of
$K$, and, if $\Lambda$ is fixed, the $2$-fold bimodule automorphisms
of $\Lambda$.  

It is important to realize that the tensor product
$\Gamma _{Y,X}\otimes X\otimes Y$ does not use the right module structures
on $X$ and $Y$, only the left module structures.  Thus these right
module structures are still available to make $\Gamma _{Y,X}\otimes
X\otimes Y$ into a $2$-fold bimodule.  

One could similarly prove that natural transformations between
additive symmetric monoidal functors represented by $X_{1}$ and
$X_{2}$ are induced by maps of bimodules $X_{1}\xrightarrow{}X_{2}$.  

\begin{proof}
Suppose $F$ is a symmetric monoidal functor from $\wedge$ to $\Box$
with a right adjoint.  Then Watts' theorem implies that there is a
bimodule $X$ and a natural isomorphism $X\otimes M\xrightarrow{}FM$.
Because $F$ is symmetric monoidal, we have a natural isomorphism 
\[
m_{M,N} \xrightarrow{} FM \Box FN \xrightarrow{} F (M\wedge N).  
\]
This translates to a natural isomorphism 
\[
m_{M,N}\mathcolon \Gamma_{Y,X} \otimes (X\otimes M) \otimes (Y\otimes
N) \xrightarrow{} X\otimes (\Lambda_{N,M}\otimes M\otimes N),
\]
where $Y$ is just another name for $X$.  Taking $M=N=R$ gives us the
desired isomorphism $m$.  The unit isomorphism $\eta \mathcolon
K'\xrightarrow{}FK$ is just the map $\eta \mathcolon
K'\xrightarrow{}X\otimes K$.  On the other hand, given $X$, $m$, and
$\eta$, we define $FM=X\otimes M$, $\eta$ in the obvious way, and
$m_{M,N}$ by naturality from $m$.  We leave to the reader the
translation between the compatibility diagrams of $F$ and the diagrams
in the theorem.  
\end{proof}

\section{Examples}\label{sec-basic-examples}

In this section, we consider some examples of additive closed
symmetric monoidal categories on $R$-modules.  In particular, we find
rings $R$ where there are no such structures, where there is eactly
one (up to additive symmetric monoidal equivalence), where there are
exactly seven, and where there are a proper class.

The most obvious case is when $R$ is a commutative ring, where
$-\wedge -$ is the usual tensor product.  This corresponds to
$\Lambda_{B,A}=R$ with $a\odot_{1}r=a\odot_{2}r=ar=ra$.  The maps $a$
and $c$ are both identity maps, the unit $K$ is $R$ itself, and the
map $\ell$ is multiplication.

Now suppose $R$ is a cocommutative Hopf algebra over a field $K$, with
diagonal $\Delta$ and counit $\epsilon$.  As is well-known, the
category of $R$-modules then becomes a closed symmetric monoidal
category under the functor $-\otimes_{K}-$, where $R$ acts by taking
the diagonal and having it then act on each factor.  This corresponds
$\Lambda_{B,A}=R_{B}\otimes_{K}R_{A}$.  The right $R$-module
structures are just right multiplication on the two factors, and the
left $R$-module structure is, as we mentioned above, the composite
\[
R\otimes_{K} (R\otimes_{K}R) \xrightarrow{\Delta \otimes 1}
R\otimes_{K} R \otimes_{K} R \otimes_{K} R \xrightarrow{1\otimes
T\otimes 1} R\otimes_{K}R\otimes_{K}R\otimes_{K}R \xrightarrow{\mu
\otimes \mu} R\otimes_{K} R,
\]
where $\Delta$ denotes the diagonal and $\mu$ denotes the
multiplication.  The commutativity isomorphism is just the twist map
$c\mathcolon
R_{B}\otimes_{K}R_{A}\xrightarrow{}R_{A}\otimes_{K}R_{B}$.  For this
to be a map of left $R$-modules, we need $R$ to be cocommutative.  The
unit isomorphism is the obvious isomorphism
\[
\ell \mathcolon (R_{B}\otimes_{K}R_{K})\otimes_{R}K\cong
R_{B}\otimes_{K} K \cong R_{B}.  
\]
The associativity isomorphism 
\[
a\mathcolon (R_{C}\otimes_{K} R_{R})\otimes_{R} (R_{B}\otimes_{K}
R_{A}) \xrightarrow{} (R_{R}\otimes_{K}R_{A})\otimes_{R}
(R_{C}\otimes_{K}R_{B})
\]
is more confusing.  It is pretty clear that we should define 
\[
a (x\otimes 1\otimes z\otimes w) = 1\otimes w\otimes x\otimes z. 
\]
But then this forces us to define 
\[
a (x\otimes y\otimes z\otimes w) = \sum a (x\otimes 1\otimes
y'z\otimes y''w) = \sum 1\otimes y''w \otimes x\otimes y'z.   
\]
Coassociativity then implies of the diagonal on $R$ then implies,
after some painful checking, that this is a map of left $R$-modules.
We leave to the excessively diligent reader the check that all the
required coherence diagrams commute.  

We would now like to classify all the additive closed symmetric
monoidal structures on $R$-modules, up to additive symmetric monoidal
equivalence, for various $R$.  The easiest case is when the unit of
the symmetric monoidal structure is $R$ itself.  This forces $R$ to be
commutative, and in this case there is only one such closed symmetric
monoidal structure.

\begin{proposition}\label{prop-commutative}
Suppose $R$ is a ring equipped with an additive closed symmetric
monoidal structure on the category of $R$-modules with unit isomorphic
to $R$.  Then $R$ is commutative and this closed symmetric monoidal
structure is additively symmetric monoidal equivalent to the usual
one.
\end{proposition}

We point out as a general rule that if the unit $K'$ is isomorphic to
some $R$-module $K$, we can always assume that the unit is $K$, up to
symmetric monoidal equivalence.  Indeed, we construct a new additive
closed symmetric monoidal structure by leaving everything the same
except the unit isomorphism $\ell $, which we modify by the
isomorphism so that the unit is $K$.  The coherence diagrams still
commute, so this is a closed symmetric monoidal structure.  A
symmetric monoidal equivalence between this new structure and the old
one is given by the identity functor, with the unit map
$K'\xrightarrow{}K$ given by the isomorphism.  

This proposition is saying that the coherence isomorphisms $a$,
$\ell$, and $c$ of the standard symmetric monoidal structure are
determined up to symmetric monoidal equivalence.  In fact, the proof
shows that $a$ and $c$ are exactly determined, though there is some
room for flexibility in $\ell$.

There is a quick proof that $R$ must be commutative, since the
endomorphism ring of the unit in a symmetric monoidal category must
always be commutative, and if the unit is isomorphic to $R$ that
endomorphism ring is $R$ as well.  However, this also falls out of the
coherence isomorphisms, so we re-prove this fact in the proof below.

\begin{proof}
The unit isomorphism shows that $\Lambda_{1,2}\cong R$ as a bimodule,
using the first right module structure on $\Lambda$.  We can then
assume it is $R$ using a symmetric monoidal equivalence.  Define
$\sigma \mathcolon R\xrightarrow{}R$ by $\sigma (x)=1\odot_{2}x$.
Note that
\[
z\odot_{2}x= (z\cdot 1)\odot_{2}x=z (1\odot_{2}x) z\sigma (x),
\]
so that $\sigma$ gives us complete information on the bimodule
$\Lambda_{1,2}$.  A similar computation shows that
$\sigma$ is a ring homomorphism.  

To see that $\sigma$ is in fact an isomorphism, consider the
commutativity isomorphism $c\mathcolon
\Lambda_{A,B}\xrightarrow{}\Lambda_{B,A}$.  This has the property that
\[
c (y)=c (1\odot_{1}y)=c (1)\odot_{2}y=c (1)\sigma (y).
\]
Since $c$ is an isomorphism, we conclude that $\sigma$ is an
isomorphism.  

We claim that associativity forces $R$ to be commutative and $\sigma$
to be the identity.  Indeed, we have
\[
a\mathcolon \Lambda_{C, \Lambda} \otimes \Lambda_{B,A} \xrightarrow{}
\Lambda_{\Lambda ,A} \otimes \Lambda_{C,B}.  
\]
Both of these are isomorphic to $R$ as left modules.  In the domain, we have 
\[
x\otimes y=x (1\otimes y)=x (\sigma (y)\otimes 1) = x\sigma (y) (1\otimes 1),
\]
and in the target we have 
\[
z\otimes w = z (1\otimes w)=z (w\otimes 1)= zw (1\otimes 1).  
\]
Thus $a$ is determined by $a (1\otimes 1)$, which must be $\gamma
\otimes 1$ for some unit $\gamma \in R$.  We will then have 
\[
a (x\otimes y)= x\sigma (y)(\gamma \otimes 1).  
\]

We first show that $R$ is commutative.  Choose an arbitrary $r,s\in
R$.  Find $y$ such that $\sigma (y)=s\gamma^{-1}$, using the fact that
$\sigma$ is an isomorphism.  Then we have 
\[
a ((1\odot_{C}r)\otimes y) = a (1\otimes y)\odot_{C}r.  
\]
But we have 
\[
a ((1\odot_{C}r)\otimes y) = a (r\otimes y)=r\sigma (y)\gamma
(1\otimes 1) = rs (1\otimes 1), 
\]
and 
\[
a (1\otimes y)\odot_{C}r = \sigma (y)\gamma (1\otimes r) =s (r\otimes
1)=sr (1\otimes 1).  
\]
We conclude that $rs=sr$, so $R$ is commutative.  

We must also have 
\[
a (1\otimes (1\odot_{A}z)) = a (1\otimes 1)\odot_{A}z.
\]
This means that 
\[
a (1\otimes \sigma (z)) = (\gamma \odot_{A}z) \otimes 1 \text{ so }
\sigma^{2} (z)\gamma (1\otimes 1) = \gamma \sigma (z) (1\otimes 1).  
\]
Thus $\sigma^{2} (z)=\sigma (z)$ for all $z$.  Since $\sigma$ is
necessarily one-to-one, we conclude that $\sigma (z)=z$. 

We now know that $\Lambda_{1,2}$ is isomorphic to $R$ with both right
module structures, and the left module structure, equal to the
canonical one.  Then the associativity isomorphism $a$ of
Theorem~\ref{thm-coherence} is just an isomorphism of left $R$-modules
from $R$ to itself, so must be right multiplication by some unit $r$.
But then the associativity pentagon shows that $r^{2}=r^{3}$, so
$r=1$.  Similarly, the commutativity isomorphism is right
multiplication by a unit $s$, and, since we now know $a$ is the
identity, the commutativity-associativity hexagon says that $s^{2}=s$,
so $s=1$.  Finally, $\ell $ must also be multiplication by some unit
$t$, but the coherence diagrams will commute no matter what $t$ is.
However, we can define a symmetric monoidal equivalence from the usual
symmetric monoidal structure to the one with $\ell =t$ by letting $F$
be the identity functor, letting the natural isomorphism $m$ be the
usual one, and letting $\eta \mathcolon R\xrightarrow{}R$ be right
multiplication by $t$.
\end{proof}

There are some simple cases where $R$ is the only possible unit of a
closed symmetric monoidal structure on the category of $R$-modules.  

\begin{theorem}\label{thm-abelian}
Let $n$ be an integer.  There is a unique additive closed symmetric
monoidal category structure on the category of $\Z /n\Z$-modules, up
to symmetric monoidal equivalence.
\end{theorem}

This theorem was proved in case $n=0$ by Foltz, Lair, and
Kelly~\cite{foltz-lair-kelly}.   

\begin{proof}
A right or left $\Z/n\Z $-module structure on an abelian group is
unique; we must have $nx=xn=x+x+ \dotsb +x$ for $n\geq 0$ and the
negative of this for $n<0$.  Thus the $2$-fold bimodule $\Lambda $
needed to define a closed symmetric monoidal structure on $\Z
/n\Z$-modules is simply a $\Z /n\Z$-module, with all of the module
structures being the same.  The unit isomorphism guarantees that
$\Lambda$ is in the Picard group of $\Z /n\Z $, which is trivial
(see~\cite[Example~2.22D]{lam}).  Hence there is an isomorphism
$f\mathcolon \Lambda \xrightarrow{}\Z /n\Z$.
Proposition~\ref{prop-commutative} completes the proof.
\end{proof}

The other simple case is when $R$ is a division ring.  

\begin{theorem}\label{thm-field}
Suppose $k$ is a division ring.  If $k$ is not a field, then there is
no additive closed symmetric monoidal structure on the category of
$k$-modules.  If $k$ is a field, there is a unique additive closed
symmetric monoidal category structure on the category of $k$-modules,
up to symmetric monoidal equivalence.
\end{theorem}

\begin{proof}
Suppose we have a closed symmetric monoidal structure induced by
$\Lambda_{B,A}$.  The unit isomorphism 
\[
\Lambda_{B,K}\otimes_{k}K\cong k
\]
shows that $K$ has to be a one-dimensional vector space, so is
isomorphic to $k$.  Proposition~\ref{prop-commutative} completes the
proof.
\end{proof}

Since the axioms for an additive closed symmetric monoidal structure on the
category of $R$-modules do not actually mention $R$ itself, the
existence and number of such structures are both Morita invariant.
Hence we get the following corollary. 

\begin{corollary}\label{cor-simple-artinian}
Suppose $R$ is a simple artinian ring, so that $R\cong M_{n} (D)$ for
some division ring $D$ and some integer $n$.  If $D$ is commutative,
there is a unique additive closed symmetric monoidal structure on the
category of $R$-modules, up to symmetric monoidal equivalence.  If $D$
is not commutative, then there is no additive closed symmetric
monoidal structure on the category of $R$-modules.
\end{corollary}

The unit of the closed symmetric monoidal structure on $M_{n}
(k)$-modules, for $k$ a field, is the unique simple left $M_{n}
(k)$-module $k^{n}$.  

To find a case where the additive closed symmetric monoidal structure
is not unique, we consider the group ring $k[\Z /2]$.  Even in this
simple case, the classification of additive closed symmetric monoidal
structures is quite involved, and will take the rest of this section
and many lemmas.  The ring $k[\Z /2]$ is both a commutative ring and a
Hopf algebra, so we know there are at least two closed symmetric
monoidal structures.  The behavior of this group ring depends on
whether the characteristic of $k$ is $2$, so we begin with this case.

We start by identifying the Hopf algebra structures on $k[\Z /2]$.  

\begin{lemma}\label{lem-Hopf-2}
Suppose $k$ is a field of characteristic $2$, and $R=k[\Z /2]\cong
k[x]/ (x^{2})$.  There are two different isomorphism classes of Hopf
algebra structures on $k$, one represented by $H_{0}$, in which
$\Delta (x)=1\otimes x+x\otimes 1$, and one represented by $H_{1}$, in
which $\Delta (x)=1\otimes x+x\otimes 1+ x\otimes x$.  
\end{lemma}

We only use the associativity and unit axioms to prove this lemma, so
it follows that every bialgebra structure on $k[\Z /2]$ is a
cocommutative Hopf algebra structure.  

\begin{proof}
The counit $\epsilon$ of a Hopf algebra structure must have $\epsilon
(1)=1$ since it is a $k$-algebra map, and $\epsilon (x)=0$ since $x$
is nilpotent.  We must have 
\[
\Delta (x) = a_{1} (1\otimes 1) + a_{2} (1\otimes x) + a_{3} (x\otimes
1) + a_{4} (x\otimes x)
\]
for some $a_{1},a_{2},a_{3},a_{4}\in k$.  The fact that 
\[
0 = \Delta (x^{2}) = \Delta (x)^{2}
\]
implies that $a_{1}=0$.  The fact that $\Delta$ is counital implies
that $a_{2}=a_{3}=1$.  One then checks that $\Delta$ is coassociative
no matter what $a_{4}$ is.  It is of course also cocommutative, and $c
(x)=x$ defines the only possible conjugation on $k[\Z /2]$.  

Ler $R_{a}$ denote the Hopf algebra where the coefficient of $x\otimes
x$ in $\Delta (x)$ is $a$.  Any isomorphism $f\mathcolon
R_{a}\xrightarrow{}R_{b}$ of Hopf algebras must be compatible with the 
counit, from which we conclude that $f (x)=rx$ for some nonzero $r\in
k$.  But then compatibility with $\Delta$ will hold if and only if
$ra=b$.  So if $b$ and $a$ are both nonzero, $r=a/b$ will yield the
desired isomorphism, but $R_{0}$ is not isomorphic to any other
$R_{b}$.  
\end{proof}

In both of these two Hopf algebra structures on $R=k[\Z /2]$ (where
$k$ has characteristic $2$), the corresponding symmetric monoidal
structure has $\Lambda =R\otimes_{k}R$, freely generated as a left
$R$-module by $m=1\otimes 1$ and $m\odot_{2}x=1\otimes x$.  We also
have $m\odot_{1}x=x\otimes 1$.  However, in $H_{0}$ we have
\[
xm= 1\otimes x + x\otimes 1 \text{ so } m\odot_{1} x = xm + m\odot_{2}x.
\]
In $H_{1}$, though, we have 
\[
xm = 1\otimes x + x\otimes 1 + x\otimes x \text{ so } m\odot_{1} x =
xm + (1+x)m\odot_{2} x.  
\]
In both cases, the unit isomorphism $\ell \mathcolon \Lambda \otimes
k\xrightarrow{}R$ has $\ell (m)=1$ and $\ell (m\odot_{2}x)=0$.  Also
the commutativity isomorphism is defined by $c (m)=m$ (and thus $c
(m\odot_{2}x)=m\odot_{1}x$).  The associativity isomorphism has 
$a (m\otimes m) = m\otimes m$ in both cases, but in $H_{0}$ we have 
\[
a (m\odot_{2}x)\otimes m) = 1\otimes x\otimes 1\otimes 1 + 1\otimes
1\otimes 1\otimes x = m\odot_{2}x \otimes m + m\odot_{2}x \otimes m,
\]
whereas in $H_{1}$ we have 
\begin{multline*}
a (m\odot_{2}\otimes m) = 1\otimes x\otimes 1\otimes 1 + 1\otimes
1\otimes 1 \otimes x + 1\otimes x\otimes 1\otimes x  \\
= m\odot_{2}\otimes m + m\otimes m\odot_{2}x + m\odot_{2}x \otimes
m\odot_{2} x.
\end{multline*}

Let us refer to these $k$-linear closed symmetric monoidal structures as
$\wedge_{H_{0}}$ and $\wedge_{H_{1}}$.  Note that $X\wedge
Y=X\otimes_{k}Y$ as $k$-modules in either case, it is just the action
of $\Z /2$ differs.  

\begin{theorem}\label{thm-2-2}
Suppose $k$ is a field of characteristic $2$, and let $R=k[\Z /2]\cong
k[x]/(x^{2})$.  Suppose $-\wedge -$ is a $k$-linear closed symmetric
monoidal structure on the category of $R$-modules with unit $K$.  Then
one of the following must hold. 
\begin{enumerate}
\item $K\cong R$ and $-\wedge -$ is $k$-linearly equivalent to
$-\otimes -$.  
\item $K\cong k$ and $-\wedge -$ is $k$-linearly equivalent to
$-\wedge_{H_{1}}-$ as a monoidal functor, but not necessarily as a
symmetric monoidal functor.  
\item $K\otimes k$ and $-\wedge -$ is $k$-linearly equivalent to
$-\wedge_{H_{0}}-$ as a unital functor, but not necessarily as a
monoidal functor.  
\end{enumerate}
In addition, we have
\begin{enumerate}
\item The isomorphism classes $(-\wedge_{H_{1}}-,\beta)$ of closed
symmetric monoidal structures with underlying monoidal functor
$-\wedge_{H_{1}}-$ are parametrized by elements $\beta \in k$, where
$c (m)=m+\beta x (m\odot_{2}x)$ in the symmetric monoidal structure
corresponding to $\beta$.
\item The isomorphism classes $(-\wedge_{H_{0}}-,\gamma)$ of closed
monoidal structures with underlying unital functor $-\wedge_{H_{0}}-$
are parametrized by elements
\[
\gamma \in \{0 \} \cup k^{\times}/ (k^{\times})^{3},
\]
where 
\[
a (m\otimes m) = m\otimes m + \gamma x (m\odot_{2}x\otimes
m\odot_{2}x)
\]
in the monoidal structure corresponding to $\gamma$.  
\item The isomorphism classes $(-\wedge_{H_{0}}-,\gamma ,\beta)$ of
closed symmetric monoidal structures with
underlying monoidal functor $(-\wedge_{H_{0}}-, \gamma)$ are
parametrized by elements $\beta $, where 
\[
c (m) = m +\beta x (m\odot_{2}x)
\] 
in the symmetric monoidal structure corresponding to $\beta$, as
follows.  
\begin{enumerate}
\item If $\gamma =0$, then $\beta \in \{0 \}\cap k^{\times }/
(k^{\times})^{2}$.  
\item If $\gamma \neq 0$ and $k$ does not have a primitive cube root
of $1$, then $\beta \in k$.  
\item If $\gamma \neq 0$ and $k$ does have a primitive cube root
$\omega$ of $1$, then $\beta =\{0 \}$ or a coset of the action of $\Z
/3$ on $k^{\times}$ given by the action of $\omega$.
\end{enumerate}
\end{enumerate}
\end{theorem}

Just so we have a specific concrete example, this theorem says that
when $k=\Z /2$, there are seven $k$-linear isomorphism classes of
$k$-linear closed symmetric monoidal structures on $k[\Z /2]$-modules,
one corresponding to the usual tensor product, two corresponding to
different symmetric monoidal structures on $-\wedge_{H_{1}}-$, and
four corresponding to different structures on the underlying unital
functor $-\wedge_{H_{0}}-$.

We will prove this theorem through a series of lemmas. 

\begin{lemma}\label{lem-2-unit}
Let $R=k[x]/ (x^{2})$ where $k$ is a field, and suppose $-\wedge -$ is
a $k$-linear closed symmetric monoidal structure on the category of
$R$-modules, with unit $K$.  Then either $K\cong R$ or $K\cong k$.  If
$K\cong R$, then $-\wedge -$ is equivalent to the usual tensor
product.  
\end{lemma}

\begin{proof}
Every $R$-module is equivalent to a direct sum of copies of $k$ and
$R$.  Any decomposition of $K$ as a direct sum of $R$-modules induces
a decomposition of $R$ as a direct sum of $R$-bimodules, via the unit
isomorphism $\Lambda_{R,K}\otimes K\cong R_{R}$.  Since $R$ is
indecomposable, $K$ must also be indecomposable, so either $K\cong k$ or
$K\cong R$.  The last statement follows from
Proposition~\ref{prop-commutative}.  
\end{proof}

\begin{lemma}\label{lem-2-basic}
Let $R=k[x]/ (x^{2})$ where $k$ is a field, and suppose $-\wedge -$ is
a $k$-linear closed symmetric monoidal structure on the category of
$R$-modules, with unit $k$.  Let $\Lambda$ be the $2$-fold bimodule
inducing $-\wedge -$.  Then there is an element $m\in \Lambda$ such
that
\[
\Lambda \cong Rm \oplus R (m\odot_{2}x) \cong R\oplus R (m\odot_{2}x)
\]
as left $R$-modules, where $\ell (m\otimes 1)=1$ and $\ell
(m\odot_{2}\otimes 1)=0$.  Furthermore, $m\odot_{2}x\neq 0$.  
\end{lemma}

This lemma also says that $\Lambda$ is a principal bimodule under the
right action $\odot_{2}$, generated by $m$.  We do not know yet
whether $\Lambda$ is a free bimodule on $m$, though we will prove this
later.

\begin{proof}
Since $\Lambda \otimes_{B,k}\otimes k\cong R$, we have 
\[
\Lambda / (\Lambda \odot_{2}x)\cong R.  
\]
If we choose an $m\in \Lambda$ with $\ell (m\otimes 1)=1$, then for any
$\lambda \in \Lambda$, we have 
\begin{multline*}
\lambda =\ell (\lambda\otimes 1)m+ n\odot_{2}x = \ell (\lambda \otimes
1)m + (\ell (n\otimes 1)m+d\odot_{2} x) \odot_{2} x \\
= \ell (\lambda
\otimes 1) m + \ell (n\otimes 1) (m\odot_{2}x),
\end{multline*}
where $n$ and $d$ denote unknown elements of $\Lambda$.  
Thus $\Lambda$ is generated as a left $R$-module by $m$ and
$m\odot_{2}x$.  Note also that 
\[
\lambda \odot_{2} x = \ell (\lambda \otimes 1) (m\odot_{2}x).  
\]
In particular, $m\odot_{2} x\neq 0$, since if it were then $\Lambda
\odot_{2} x=0$.  The commutativity isomorphism then implies $\Lambda
\odot_{1}x=0$, so 
\[
Rx = (\Lambda \odot_{1}x) \otimes k=0,
\]
which is a contradiction.  
\end{proof}

At this point, we have not determined whether $m\odot_{2}x$ generates
a copy of $k$ or a copy of $R$.  This depends on whether $x
(m\odot_{2}x)=0$ or not.  Note, however, that $\ell
(m\odot_{1}x\otimes 1)=x$, so we must have
\[
m\odot_{1} x = xm + b (m\odot_{2}x).  
\]
for some $b\in R$.  We will then have 
\[
(m\odot_{2} x) \odot_{1} x = (xm + b (m\odot_{2}x))\odot_{2} x
=xm\odot_{2}x.  
\]
It will be helpful in what follows if we write $a=a^{0}+a^{1}x$ for
elements $a\in R$, with $a^{0},a^{1}\in k$.  

\begin{lemma}\label{lem-2-commute}
Suppose the characteristic of $k$ is $2$.  In the situation of
Lemma~\ref{lem-2-basic} and with the notation above, there is a
commutativity involution $c\mathcolon
\Lambda_{A,B}\xrightarrow{}\Lambda_{B,A}$ if and only if $b^{0}=1$,
and in that case we have
\[
c (m) = (1+\alpha x)m + (\beta x) (m\odot_{2}x)
\]
for some $\alpha ,\beta \in k$, which can be anything.  This implies 
\[
c (m\odot_{2}x) = xm + (b+\alpha x) (m\odot_{2}x).
\]  
\end{lemma}


\begin{proof}
Write
\[
c (m)=rm+ s(m\odot_{2}x),
\]
for some $r,s \in R$.  Then 
\[
c (m\odot_{2}x)=c (m)\odot_{1}x = rxm +rb
(m\odot_{2}x)+ sx(m\odot_{2}x).  
\]
In order for $c$ to have the desired properties, we need 
\[
c (m\odot_{1}x) = c (m) \odot_{2} x \text{ and } c^{2} (m)=m.  
\]
In order for $c^{2} (m)=m$, computation shows that we need 
\[
r^{2}+rsx=1 \text{ and }  (rs+rsb+s^{2}x) (m\odot_{2}x)=0.
\]
The first equation is equivalent to $(r^{0})^{2}=1$ and
$2r^{0}r^{1}+r^{0}s^{0}=0$, so $r^{0}=\pm 1$ and $s^{0}=-2r^{1}$.
Since we are in characteristic $2$, this means $r^{0}=1$ and
$s^{0}=0$.  We leave the second equation aside for the moment. In
order for $c (m\odot_{1}x)=c (m)\odot_{2}x$, computation shows that we
need
\[
r (1+b)x=0 \text{ and } (sx+rb^{2}+sbx) (m\odot_{2}x)=r
(m\odot_{2}x).  
\]
Since we know that $r^{0}=\pm 1$, this first equation implies that
$b^{0}=-1$; since we are in characterstic $2$, $b^{0}=1$.  Computation
then shows that, in characteristic $2$, these conditions guarantee
that both equations involving $m\odot_{2}x$ hold, whether $x
(m\odot_{2}x)=0$ or not.
\end{proof}

In fact, we can make $\alpha =0$ by modifying our choice of $m$.  

\begin{lemma}\label{lem-2-modify}
In the situation of Lemma~\ref{lem-2-commute}, so in particular when
the characteristic of $k$ is $2$, we can modify our choice of the
bimodule generator $m$ of $\Lambda$ so that 
\[
c (m) = m + \beta x (m\odot_{2}x)
\]
for some $\beta \in k$.  In this case, $c (m\odot_{2}x)=m\odot_{1}x$.
\end{lemma}

From now on, we assume $m$ is chosen as in Lemma~\ref{lem-2-modify}.  

\begin{proof}
Let $n=m+\alpha (m\odot_{2}x)$.  Then $\ell (n\otimes 1)=1$, so $n$ is
a perfectly good bimodule generator for $\Lambda$.  Note that 
\[
n\odot_{2}x = m\odot_{2}x \text{ and } n\odot_{1}x = xn + b (n\odot_{2}x)
\]
as before (after some calculation).  However, we have 
\[
c (n) = n + \gamma x (n\odot_{2}x),
\]
after some calculation, for some $\gamma \in k$.  
\end{proof}

We must now come to grips with the associativity isomorphism
\[
a\mathcolon \Lambda_{C, \Lambda} \otimes \Lambda_{B,A} \xrightarrow{}
\Lambda_{\Lambda , A} \otimes  \Lambda_{C,B}.
\]
If $\Lambda$ is a free bimodule generated by $m$, then both the domain
and range of $a$ are isomorphic to $R^{\oplus 4}$ as left $R$-modules,
with summands generated by $m\otimes m$, $m\odot_{2}x\otimes m$,
$m\otimes m\odot_{2}x$, and $m\odot_{2}x\otimes m\odot_{2}x$.  If
$\Lambda$ has dimension $3$, on the other hand, both the domain and
range of $a$ are isomorphic to $R\oplus k\oplus R$, with summands
generated by $m\otimes m, m\odot_{2}x \otimes m$ and $m\otimes
m\odot_{2}x$, respectively.  For example, in the domain we have
\[
m\otimes xm\odot_{2}x =m\odot_{2} x \otimes m\odot_{2}x,
\]
but in the range we have 
\[
m\otimes xm\odot_{2}x = m\odot_{1}x\otimes m\odot_{2}x = xm\otimes
m\odot_{2}x + (1+b^{1}x) m\odot_{2} x \otimes m\odot_{2} x. 
\]
With respect to this basis, write 
\[
a (m\otimes m) = (e_{1},e_{2}, e_{3}, e_{4}) \text{ and } a
(m\odot_{2}x\otimes m) = (f_{1},f_{2}, f_{3}, f_{4}).  
\]
If the dimension of $\Lambda$ is $3$, then we take $e_{4}=f_{4}=0$.
In addition, in that case $f_{3}^{0}=0$, since $x (m\otimes
m\odot_{2}x)=0$.  

\begin{lemma}\label{lem-2-unit-assoc}
With the above definitions, if $a$ is a map of $3$-fold bimodules,
then the left and right unit coherence diagram commutes if and only if
$e_{1}=1$, $e_{2}=0$, $f_{1}=0$, and $f_{2}=1$.
\end{lemma}

\begin{proof}
Apply the coherence diagram to $m\otimes m\otimes 1$ and to
$m\odot_{2}x\otimes m\otimes 1$.  
\end{proof}

We still have to determine what conditions are necessary for $a$ to be
a map of $3$-fold bimodules.  

\begin{lemma}\label{lem-2-assoc}
In order for the map $a$ defined above to be a map of $3$-fold
bimodules making the left and right unit coherence diagram commute,
$\Lambda$ must be the free bimodule on $m$, and 
\[
e_{1}=1, e_{2}=0, e_{3}^{0}=0, f_{1}=0, f_{2}=1, f_{3}=1, f_{4}=
b^{1}+ e_{3}.  
\]
\end{lemma}

\begin{proof}
Of course, we know already that $e_{1}=1, e_{2}=0, f_{1}=0,f_{2}=1$,
with $f_{3}^{0}=0$ if $\Lambda$ has dimension $3$.  We have implicitly
assumed that $a$ is a map of left $R$-modules by defining it only in
terms of generators.  To ensure that $a$ preserves the right module
structure represented by $A$ in 
\[
a\mathcolon \Lambda_{C,\Lambda} \otimes \Lambda_{B,A} \xrightarrow{}
\Lambda_{\Lambda ,A} \otimes \Lambda_{C,B}, 
\]
we must have
\[
a (m\otimes m\odot_{2}x) = (1,0,e_{3},e_{4})\odot_{A} x = (0,1,0,e_{3})
\]
and 
\[
a (m\odot_{2}x\otimes m\odot_{2}x) = (0,0,0,f_{3}).  
\]
We now turn to the right module structure represented by $B$.  Here we
must have
\[
a (m\otimes m\odot_{1}x)=a (m\otimes m)\odot_{B}x.
\]
Calculation shows that this forces $f_{3}=1$, and this rules out the
case when the dimension of $\Lambda$ is $3$ (since $f_{3}^{0}=0$ if
the dimension of $\Lambda$ is $3$).  The same calculation shows that
$f_{4}=b^{1}+e_{3}$.  Further calculation shows that this is enough to
ensure that $a$ preserves the right module structure represented by
$B$.  We now must ensure that $a$ preserves the right module structure
represented by $C$.  More calculation of the relation
\[
a (m\odot_{1}x\otimes m)=a (m\otimes m)\odot_{C}x
\]
gives $e_{3}x=0$, so $e_{3}^{0}=0$.  Further calculation implies that
this enough to make $a$ preserve the right module structure
represented by $C$.
\end{proof}

We turn finally to the associativity pentagon.  

\begin{lemma}\label{lem-2-assoc-pentagon}
Given that $a$ satisfies the conditions of Lemma~\ref{lem-2-assoc},
$a$ makes the associativity pentagon commute if and only if 
\[
e_{3}=0, f_{4}=b^{1}, e_{4}^{0}=0 \text{ and either } e_{4}^{1}=0
\text{ or } b^{1}=0.  
\]
\end{lemma}

\begin{proof}
If we apply the associativity pentagon to $m\odot_{2}x\otimes m\otimes
m$, we eventually find
\[
e_{3}=0 \text{ and } e_{4}^{0} = 0.
\]
and so $f_{4}=b^{1}+e_{3}=b^{1}$.  Further computation with the
associativity pentagon applied to $m\otimes m\otimes m$ eventually
yields
\[
e_{4}^{1} b^{1}=0, \text{ so } e_{4}^{1}=0 \text{ or } b^{1}=0.
\]
These conditions then make the associativity pentagon commute.  
\end{proof}

We are now left with determining which of these different $\Lambda$
define additively equivalent symmetric monoidal structures.  Recall
that for such an equivalence, we need an element $X$ of the bimodule
Picard group, an isomorphism $\eta \mathcolon X\otimes
k\xrightarrow{}k$, and an isomorphism
\[
\Gamma \otimes X\otimes X \xrightarrow{} X\otimes \Lambda 
\]
making various diagrams commute.  

\begin{lemma}\label{lem-2-bimodule-Picard}
Let $R=k[x]/x^{2}$.  Up to isomorphism, the only invertible
$R$-bimodules are the $R_{u}$, where $u$ is a unit in $k$, $R_{u}=R$
as a left module, and $1\odot x=ux$.
\end{lemma}

This follows immediately from~\cite[Lemma~3.3]{yekutieli}.  With all
these lemmas in hand, we can now complete the proof of
Theorem~\ref{thm-2-2}.  

\begin{proof}[Proof of Theorem~\ref{thm-2-2}]
Let $X=R_{u}$, and suppose we have an isomorphism $\eta \mathcolon
k\xrightarrow{}X\otimes k$ of left modules and an isomorphism
\[
q \mathcolon \Gamma \otimes X \otimes X\xrightarrow{}X\otimes \Lambda 
\]
of $2$-fold bimodules making the compatibility diagrams of
Theorem~\ref{thm-equiv} commute.  The map $\eta$ is determined by
$\eta (1)=1\otimes \rho $ for some nonzero $\rho$ in $k$, and the map
$q$ is determined by $q (m\otimes 1\otimes 1)=\sigma (1\otimes m)+\tau
(1\otimes m\odot_{2}x)$ with $\sigma ,\tau \in R$, where we have used
$m$ for a bimodule generator in both $\Gamma$ and $\Lambda$ satisfying
the condition on $c (m)$ as in Lemma~\ref{lem-2-modify}. Since $q$ is
a map of bimodules, it follows that
\[
q (m\odot_{2}x\otimes 1\otimes 1) = u^{-1}q (m\otimes 1\cdot x \otimes
1) = u^{-1}\sigma (1\otimes m\odot_{2}x).
\]
In particular, for $q$ to be an isomorphism, we need $\sigma$ to be a
unit in $R$.  We also need 
\[
q (m\odot_{1}x \otimes 1\otimes 1) = u^{-1}q (m\otimes 1\otimes 1\cdot
x) = u^{-1}q (m\otimes 1\otimes 1)\odot_{1}x.
\]
Further calculation with this last equation yields
$b^{1}_{\Lambda}=u^{-1}b^{1}_{\Gamma}$.  

This means that if $b^{1}_{\Gamma }$ is nonzero, we can choose $u$ so
as to ensure $b^{1}_{\Lambda}=1$.  Said another way, given $\Gamma$
with $b^{1}$ nonzero, we can define $\Lambda =R_{u^{-1}}\otimes \Gamma
\otimes R_{u}\otimes R_{u}$ for a suitable $u$ and get an additively
equivalent symmetric monoidal structure with $b^{1}=1$.  Any
isomorphism between a $\Gamma$ and a $\Lambda$ both with $b^{1}=1$
must have $u=1$, so we can think of it as an automorphism of $\Lambda
$.  It is still useful to use $\Gamma$ for the domain copy of $\Lambda
$, because the choice of generator $m$ as in Lemma~\ref{lem-2-modify}
could be different in the two copies of $\Lambda$.  We can then work
through the compatibility diagrams of Theorem~\ref{thm-equiv} in this
case.  The unit diagram forces $\sigma =\rho^{1}$, so $\sigma \in k$.
In order for the commutativity diagram to commute, we need $\tau x=0$,
so $\tau^{0}=0$.  Then one finds that the $\beta$ in the commutativity
isomorphism in both $\Gamma$ and $\Lambda$ must be the same.  Since we
are in the case where $b^{1}=1$, the associativity isomorphism is
completely determined by th preceding lemmas.  Thus we find that if
$b^{1}$ is nonzero, then the monoidal structure determined
by $\Lambda$ is additively equivalent to the one given by the Hopf
algebra $H_{1}$.  Since $\beta$ does not change under these
isomorphisms, it could be anything, so the different symmetric
monoidal structures on $-\wedge_{H_{1}}-$ are classified by $\beta$.  

Now suppose $b^{1}_{\Gamma}=b^{1}_{\Lambda}=0$.  As above, we must
have $\sigma =\rho^{-1}$ in order to ensure the unit compatibility
diagram commutes.  The commutativity compatibility diagram again
forces $\tau^{0}=0$, but this time we have
$\beta_{\Lambda}=u^{-2}\beta_{\Gamma}$.  Since $b^{1}=0$, the
preceding lemmas allow for a nontrivial $e_{4}^{1}$ as well, so we
must check the associativity compatibility diagram too.  Painful
computation then gives that the $e_{4}^{1}$ for $\Lambda$ is $u^{-3}$
times the $e_{4}^{1}$ for $\Gamma$.  Hence the different monoidal
structures on $-\wedge_{H_{0}}-$ are classified by 
\[
\gamma \in \{0 \} \cup k^{\times}/ (k^{\times})^{3}.  
\]
If $\gamma =0$, then we can use $R_{u}$ to make the symmetric monoidal
structure corresponding to $\beta $ isomorphic to the one
corresponding to $u^{-2}\beta$.  Hence the symmetric monoidal
structures when $\gamma =0$ correspond to 
\[
\beta \in \{0 \}\cup k^{\times}/ (k^{\times})^{2}.
\]
When $\gamma \neq 0$, if we fix $\gamma$ we can only use the $R_{u}$
with $u^{3}=1$.  If $k$ has no primitive cube root of $1$, then, the
symmetric monoidal structures are parametrized by $\beta \in k$, but
if $k$ does have a primitive cube root of $1$, the symmetric monoidal
structures are parametrized by orbits of $\Z /3$ acting on $k$ by
multiplying by the primitive cube root of $1$. 
\end{proof}

We now consider closed symmetric monoidal structures on $k[\Z
/2]$-modules when the characteristic of $k$ is not $2$.  Here the
answer is wildly different; there are a proper class of such
structures!

\begin{theorem}\label{thm-Z2}
Suppose $k$ is a field whose characteristic is not $2$, and let
$R=k[\Z /2]\cong k[x]/ (x^{2}-1)$.  If $-\wedge -$ is a closed
$k$-linear symmetric monoidal structure on the category of
$R$-modules, then its unit $K$ is isomorphic to $R$, $k_{+}$, or
$k_{-}$.  If the unit is isomorphic to $R$, then $-\wedge -\cong
-\otimes -$ as $k$-linear symmetric monoidal functors.  If $K\cong
k_{-}$, then $- \wedge -$ is $k$-linearly symmetric monoidal
equivalent to a closed $k$-linear symmetric monoidal structure whose
unit is $k_{+}$.  Given any $R$-module $M$, there is a $k$-linear
symmetric monoidal structure for which the unit is $k_{+}$ and
$k_{-}\wedge k_{-}=M$.
\end{theorem}

It is easy to see that symmetric monoidal structures with
nonisomorphic values of $M$ cannot be equivalent.  However, we do not
know if there is more than one closed symmetric monoidal structure for
a given $M$.

\begin{proof}
Let $x$ denote the element $[1]$ of $k[\Z /2]$, so $R=k[\Z /2]\cong
k[x]/ (x^{2}-1)$.  Since the characteristic is not $2$, this ring is
semisimple.  Any module $M$ splits as $M_{+}\oplus M_{-}$, where
$M_{+}$ is the $1$-eigenspace of $x$ and $M_{-}$ is the
$-1$-eigenspace of $x$.  In particular, $R$ itself so splits, with the
splitting given by the orthogonal idempotents $e_{+}=(1/2) (1+x)$ and
$e_{-}= (1/2) (1-x)$.  This produces a splitting 
\[
R\Mod \cong k\Mod \times k\Mod ,
\]
up to equivalence, of the entire category of $R$-modules.  Thus every
$R$-module is a direct sum of copies of $k_{+}$ and $k_{-}$, and there
are no maps from $k_{+}$ to $k_{-}$, or from $k_{-}$ to $k_{+}$.    

Given a $k$-linear closed symmetric monoidal
structure, the corresponding bimodule $\Lambda$ splits into $8$
different spaces $\Lambda_{a,b,c}$, where $a,b,c$ are each $+$ or $-$
(by simultaneously diagonalizing the $3$ actions of $x$).  Here the
$a$ stands for the left action of $R$ on $\Lambda$, and $b$ and $c$
for the two right actions.  Of course the unit $K=K_{+}\oplus K_{-}$
as well.  One can easily check that $k_{a,b,c}\otimes_{c}k_{c}\cong
k_{a,b}$, and $k_{a,b,c}\otimes_{c}k_{d}=0$ if $d\neq c$.  Since
$R=k_{+,+}\oplus k_{-,-}$, we see that there must be dimension one
terms $\Lambda_{+,+,a}$ and $K_{a}$ and dimension one terms
$\Lambda_{-,-,b}$ and $K_{b}$.

If $a\neq b$, then $K=R$, and so Proposition~\ref{prop-commutative}
tells us that our closed symmetric monoidal structure is equivalent to
$\otimes_{R}$.  We can therefore assume $a=b$.  The commutativity
isomorphism tells us $\Lambda_{x,y,z}$ has the same dimension as
$\Lambda_{x.z.y}$.  Therefore, taking $z\neq a$, we find that
$\Lambda_{z,a,z}$ is nonzero.  If $K_{z}$ were also nonzero, we would
get a term $k_{z,a}$ in $\Lambda_{B,K} \otimes K\cong R_{R}$.  Since
we cannot have such a term, we conclude that $K=k_{a}$, so the unit is
one-dimensional over $k$.  

There is an obvious self-equivalence of the category of $R$-modules
that permutes the two copies of $k\Mod$.  That is, it sends $k_{+}$ to
$k_{-}$, and vice versa.  Up to symmetric monoidal equivalence, then,
we can assume the unit of our symmetric monoidal structure is
$k_{+}$.  Let $M=k_{-}\wedge k_{-}$.  Then, using the decomposition of
$R=k_{+,+}\oplus k_{-,-}$ as bimodules, we get 
\[
\Lambda \cong k_{+,+,+} \oplus k_{-,-,+} \oplus k_{-,+,-} \oplus M_{-,-},
\]
where $M_{-,-}$ is the $2$-fold bimodule whose underlying left module
is $M$, and where $x$ acts as $-1$ in both right module structures.  

Now suppose we are given $M$.  We want to construct a $k$-linear
closed symmetric monoidal structure on $R$-modules with $k_{-}\wedge
k_{-}=M$.  We simply define $k_{+}\wedge N=N\wedge k_{+}=N$ for any
$R$-module $N$, and $k_{-}\wedge k_{-}=M$.  On morphisms, we note that
there are no nonzero morphisms from $k_{+}$ to $k_{-}$ or vice versa,
and that every endomorphism of $k_{+}$ or $k_{-}$ is given by
multiplication by an element of $k$.  So we define the induced
morphism to be multiplication by the same element of $k$.  Since every
$R$-module is a direct sum of copies of $k_{+}$ and $k_{-}$, this
defines $-\wedge -$ as a bifunctor.  We define the left unit to be the
identity.  The commutativity isomorphism 
\[
c_{xy}\xrightarrow{}k_{x}\wedge k_{y}\xrightarrow{}k_{y}\wedge k_{x}
\]
is the identity, where $x$ and $y$ denote signs.  We then extend
through direct sums.  The associativity isomorphism 
\[
a_{x,y,z}\mathcolon (k_{x}\wedge k_{y})\wedge k_{z} \xrightarrow{}
k_{x}\wedge (k_{y}\wedge k_{z}) 
\]
is the identity as long as at least one of $x$, $y$, or $z$ is $+$.
The map 
\[
a_{-,-,-}\mathcolon M\wedge k_{-} \xrightarrow{} k_{-}\wedge M
\]
is the commutativity isomorphism.  We leave to the reader the check
that the coherence diagrams hold.  
\end{proof}

\section{Structural results}\label{sec-structure}

In the examples in the previous section, especially in
Theorem~\ref{thm-Z2}, we saw that the $2$-fold bimodule $\Lambda$ can
be very complicated.  It does not have to be finitely generated, or
even countably generated.  However, it cannot be completely random
either.  Furthermore, in all the examples we have, the
unit $K$ of an additive symmetric monoidal structure on the category
of $R$-modules is always a prinicipal $R$-module.  It is tempting to
wonder whether this always holds, or whether there are other
properties that $K$ must have.  

In this section, we show that $K$ is always a finitely generated
module with commutative endomorphism ring, and that $\Lambda$ is
faithful in a very strong sense.  

We begin by noting that tensoring with $\Lambda$ reflects any property
of morphisms that the tensor product preserves.  

\begin{proposition}\label{prop-reflect}
Suppose $R$ is a ring and $\Lambda_{1,2}$ is a $2$-fold $R$-bimodule
that determines a closed symmetric monoidal structure on the category
of $R$-modules with unit $K$.  Let $\cat{P}$ be a replete class of
morphisms of abelian groups with the property that if $f\in \cat{P}$,
then $A\otimes f$ and $f\otimes B\in \cat{P}$ for any left $R$-module
$A$ or right $R$-module $B$.  If $f$ is a morphism of left
$R$-modules, then $f\in \cat{P}$ if and only if $f\otimes \Lambda \in
\cat{P}$.  Similarly, if $g$ is a morphism of right $R$-modules, then
$g\in \cat{P}$ if and only if $\Lambda \otimes g\in \cat{P}$.  This
statement holds with either right module structure on $\Lambda$.
\end{proposition}

Recall that a class of morphisms is \textbf{replete} whenever $f\in
\cat{P}$ and $f\cong g$ in the category of morphisms, then $g\in
\cat{P}$.  Said another way, if we have a commutative square 
\[
\begin{CD}
A @>f>> B \\
@ViVV  @VVjV \\
A' @>>g> B'
\end{CD}
\]
where $i,j$ are isomorphisms, then $f\in \cat{P}$ if and only if $g\in
\cat{P}$.  

\begin{proof}
Suppose $g$ is a morphism of right $R$-modules.  By definition, if
$g\in \cat{P}$, then $\Lambda \otimes g\in \cat{P}$.  Conversely,
suppose $\Lambda \otimes g\in \cat{P}$, where we use the leftmost
right module structure on $\Lambda$.  Then
\[
g\cong R_{R}\otimes g \cong (\Lambda_{B,K}\otimes K) \otimes g \cong
(\Lambda_{B,K}\otimes g)\otimes K,
\]
and this is in $\cat{P}$ since $\Lambda_{B,K}\otimes g$ is so.  We use
the commutativity isomorphism to prove the same thing for the
other right module structure on $\Lambda$.  
\end{proof}

Taking the class $\cat{P}$ to tbe the collection of zero morphisms,
the collection of isomorphisms, and the collection of surjective maps
gives us the following corollary.  

\begin{corollary}\label{cor-reflect}
Suppose $R$ is a ring and $\Lambda_{1,2}$ is a $2$-fold $R$-bimodule
that determines a closed symmetric monoidal structure on the category
of $R$-modules with unit $K$.  Let $f$ denote a morphism of right
$R$-modules and let $g$ denote a morphism of left $R$-modules.  
\begin{enumerate}
\item $f\otimes \Lambda =0$ if and only if $f=0$.  With either right
module structure, $\Lambda\otimes g=0$ if and only if $g=0$.  
\item $f\otimes \Lambda$ is an isomorphism if and only if $f$ is so.
With either right module structure, $\Lambda \otimes g$ is an
isomorphism if and only if $g$ is so.  
\item $f\otimes \Lambda$ is a surjection if and only if $f$ is so.
With either right module structure, $\Lambda \otimes g$ is a
surjection if and only if $g$ is so.
\end{enumerate}
\end{corollary}

In particular, these imply that $\Lambda$ is faithful.  

\begin{corollary}\label{cor-faithful}
Suppose $R$ is a ring and $\Lambda_{1,2}$ is a $2$-fold $R$-bimodule
that determines a closed symmetric monoidal structure on the category
of $R$-modules.  Then $\Lambda$ is faithful as a left or right
$R$-module, with either right module structure.  
\end{corollary}

\begin{proof}
Choose $r\neq 0\in R$.  Then the map $R\xrightarrow{}R$ that is left
multiplication by $r$ induces left multiplication by $r$ on $\Lambda
=R\otimes \Lambda$.  Since $r\neq 0$, this map is also nonzero by
Corollary~\ref{cor-reflect}.  Hence $r$ does not annihilate $\Lambda$,
so $\Lambda$ is faithful.  Use right multiplication by $r$ to see that
$\Lambda$ is faithful as a right $R$-module.
\end{proof}

They also imply that $K$ is finitely generated.  

\begin{theorem}\label{thm-unit-finitely-generated}
The unit $K$ in an additive closed symmetric monoidal structure on the
category of $R$-modules is finitely generated.
\end{theorem}

\begin{proof}
In the unit isomorphism $\Lambda_{B,K}\otimes K\cong R_{R}$, write $1$
as the image of a finite sum $\sum \lambda_{i}\otimes k_{i}$.  Let
$K'$ denote the submodule of $K$ generated by the $k_{i}$, and
$j\mathcolon K'\xrightarrow{}K$ denote the inclusion.  Then $\Lambda
\otimes j$ is surjective, so $j$ must also be.  
\end{proof}

We suspect that the unit $K$ must in fact be a principal $R$-module,
but we do not know how to prove this.  

Another essential property of $K$, or the unit of any symmetric
monoidal category, is that its endomorphisms commute with each other.
Somewhat more is true in our case.  

\begin{theorem}\label{thm-unit-commute}
Suppose $K$ is the unit of an additive closed symmetric monoidal
structure on the category of $R$-modules.  Then $\End_{R} (K)$ is a
subring of the center $Z (R)$ of $R$.  
\end{theorem}

\begin{proof}
Suppose $f\in \End_{R}(K)$.  Then $\Lambda \otimes f$ is a bimodule
endomorphism of $R$, through the unit isomorphism.  Any bimodule
endomorphism of $R$ must be given by $x\mapsto rx$ for some $r\in Z
(R)$.  This defines a ring homomorphism $\End_{R} (K)\xrightarrow{}Z
(R)$.  If $f$ is in the kernel of this homomorphism, then $\Lambda
\otimes f=0$, but then $f=0$ by Corollary~\ref{cor-reflect}.  
\end{proof}

We note that $\End_{R} (K)$ can be a proper submodule of $Z (R)$, as
for example when $R=k[x]/ (x^{2})$ and the unit is $k$ (of
characteristic $2$).  

It is pretty rare for an $R$-module to have a commutative endomorphism
ring.  Using the work of Vasconcelos~\cite{vasconcelos}, for example,
we can deduce the following corollary.

\begin{corollary}\label{cor-Noetherian}
Suppose $R$ is a commutative Noetherian ring with no nonzero nilpotent
elements, and $K$ is the unit of an additive closed symmetric monoidal
structure on the category of $R$-modules.  Then 
\[
K\cong \ideal{a}/\ideal{b}
\]
for some radical ideal $\ideal{b}$ and some ideal $\ideal{a}\supseteq
\ideal{b}$ with the ideal quotient $(\ideal{b}\uc
\ideal{a})=\ideal{b}$.  
\end{corollary}

Recall that the ideal quotient $(\ideal{b}\uc
\ideal{a})$ is the set of all $x$ such that $x\ideal{a}\subseteq
\ideal{b}$.  

\begin{proof}
Let $\ideal{b}=\ann (K)$, the annihilator of $K$.  We have an obvious
monomorphism of rings $R/\ideal{b}\xrightarrow{}\End_{R} (K)$ that
takes $r$ to multiplication by $r$, But $\End_{R} (K)$ is a subring of
$R$ by Theorem~\ref{thm-unit-commute}.  Hence $R/\ideal{b}$ is a
subring of $R$, and therefore has no nilpotents, so $\ideal{b}$ is a
radical ideal.  Thus $K$ is a finitely generated (by
Theorem~\ref{thm-unit-finitely-generated}), faithful
$R/\ideal{b}$-module with commutative endomorphism ring, and
$R/\ideal{b}$ is a commutative Noetherian ring with no nonzero
nilpotent elements.  Vasconcelos~\cite{vasconcelos} proves in this
situation that $K$ is an ideal in $R/\ideal{b}$.  Hence $K\cong
\ideal{a}/\ideal{b}$ for some ideal $\ideal{a}$ of $R$.  The condition
on the ideal quotient is so that the annihilator of $K$ will in fact
be $\ideal{b}$.
\end{proof}

Note that, if $M$ is a submodule of the unit $K$, then the image of
$\Lambda \otimes M$ in $\Lambda \otimes K\cong R$ will be a
sub-bimodule of $R$, and hence a two-sided ideal.  

\begin{corollary}\label{cor-simple}
Suppose $K$ is the unit of an additive closed symmetric monoidal
structure on the category of $R$-modules.  Then every nonzero proper
submodule of $K$ gives rise to a nonzero proper two-sided ideal of
$R$.  Hence, if $R$ is a simple ring, then $K$ is a simple left
module.
\end{corollary}

We note that a simple commutative ring is of course a field, but there
are many simple noncommutative rings that are not division rings.  We
would like to be able to say that the map from nonzero proper
submodules of $K$ to two-sided ideals of $R$ is one-to-one, but we do
not know if this is true.  

\begin{proof}
Suppose $L$ is a proper submodule of $K$.  Then the maps
$L\xrightarrow{}K$ and $K\xrightarrow{}K/L$ are both nonzero, so they
remain so after tensoring with $\Lambda_{B,K}$ by
Corollary~\ref{cor-reflect}.  Hence the image of
$\Lambda_{B,L}\otimes L$ is a nonzero proper subbimodule of
$\Lambda_{B,K}\otimes K=R_{R}$.  
\end{proof}

As above, we do not know if this map from submodules of $K$ to
two-sided ideals in $R$ is one-to-one, but it is on direct summands of
$K$.  

\begin{corollary}\label{cor-summands}
Suppose $K$ is the unit of an additive closed symmetric monoidal
structure on the category of $R$-modules.  There is a one-to-one map
from isomorphism classes of direct summands of $K$ to central
idempotents in $R$.  In particular, if $R$ is indecomposable as a
ring, then $K$ is indecomposable as an $R$-module.
\end{corollary}

\begin{proof}
Suppose $M$ is a direct summand of $K$, so that there is a retraction
$f\mathcolon K\xrightarrow{}M$.  Tensoring with $\Lambda$ gives us a
retraction of bimodules
$R\xrightarrow{}\Lambda \otimes M$.  The composite 
\[
R\xrightarrow{}\Lambda \otimes M\xrightarrow{}R
\]
must be multiplication by a central idempotent $e$ of $R$, with
$\Lambda \otimes M=eR$.  The bimodule $\Lambda \otimes M$ determines
$e$~\cite[Exercise~22.2]{lam-first-course}, and we can recover $M$
from $\Lambda \otimes M$, up to isomorphism, by tensoring witk $K$.  
\end{proof}


\providecommand{\bysame}{\leavevmode\hbox to3em{\hrulefill}\thinspace}
\providecommand{\MR}{\relax\ifhmode\unskip\space\fi MR }
\providecommand{\MRhref}[2]{%
  \href{http://www.ams.org/mathscinet-getitem?mr=#1}{#2}
}
\providecommand{\href}[2]{#2}

\end{document}